\documentclass[11pt]{amsart}

\usepackage[margin=1in]{geometry}
\usepackage{amssymb,mathtools}
\usepackage{enumitem}
\usepackage{hyperref}
\usepackage{microtype}

\newtheorem{theorem}{Theorem}[section]
\newtheorem{proposition}[theorem]{Proposition}
\newtheorem{lemma}[theorem]{Lemma}

\theoremstyle{definition}
\newtheorem{definition}[theorem]{Definition}
\theoremstyle{remark}
\newtheorem{remark}[theorem]{Remark}

\newcommand{\diag}{\operatorname{Diag}}
\newcommand{\charpoly}{\chi}

\title[Isospectral Deformations in Dimension Four]{Isospectral Deformations of Diagonal Matrices by Symmetric Matrices in Dimension Four}
\author{Hari Asokan}
\date{\today}
\subjclass[2020]{11C20, 11E04, 11E16, 14L24, 14L30, 15A18, 15A24, 15A22}
\keywords{isospectral deformations, symmetric matrices, diagonal matrices, characteristic polynomial}

\begin{document}

\begin{abstract}
For a diagonal matrix $D$ with distinct diagonal entries in an algebraically closed field $K$ we study the symmetric matrices $Q$ with entries in $K$ such that the affine line $D+tQ$ has the same eigenvalues in $K$ as $D$. We study the equality of characteristic polynomials, $\charpoly_{D+tQ}(x)=\charpoly_D(x) \text{  for all }\: t\in K$ and examine the resulting system of equations. For generic diagonal matrices of size at most $4$ we show that only the trivial symmetric matrix $Q=0$ satisfies this isospectral condition. We use this result to compute the degree of a field extension that appears in the invariant theory of quadratic forms under congruence action. 
\end{abstract}

\maketitle

\section{Introduction}
 Let $p\in \{0,2,3,5,7,\dots\}$. Let $K$ be an algebraically closed field of characteristic $p$; so $\operatorname{char}(K)=p.$ We assume $\operatorname{char}(K)\neq 2.$  A $K$-variety is a separated scheme of finite type over $K.$ We will denote the $K$-valued points of a $K$-variety $\mathbf{X}$ by $\mathbf{X}(K)$. By $\mathbf{M}_n$ and $\mathbf{V}_n$ we denote the $K$-variety of matrices and symmetric matrices respectively of size $n\times n$. By $\mathbf{Z}_n$ we denote the $K$-variety of diagonal matrices of size $n\times n$. By $\mathbb{P}\mathbf{V}_n$ and $\mathbb{P}\mathbf{Z}_n$ we denote the respective projective varieties. By $V_n:=\mathbf{V}_n(K)$ and $Z_n:=\mathbf{Z}_n(K)$ we denote the respective $K$-valued points, whereas $\mathbb{P}V_n:=\mathbb{P}\mathbf{V}_n(K)$ and $\mathbb{P}Z_n:=\mathbb{P}\mathbf{Z}_n(K)$ denote the $K$-valued points of the projective varieties. We denote the characteristic polynomial of the $n\times n$ matrix $A$ by $\chi_A(x).$\\
The coefficients of the characteristic polynomial  define
the quotient morphism
\[ \chi:
\mathbf{M}_n\longrightarrow
\operatorname{Spec}\!\left(K[\mathbf{M}_n]^{\operatorname{GL}_n}\right)
\simeq \mathbb A^n,
\]
where $\operatorname{GL}_n(K)$ acts by conjugation. The fibers of this
morphism are the isospectral varieties. It is of interest to study which affine sets, can lie in  a fiber. 

For a diagonal matrix
\[
D=\diag(d_1,\dots,d_n) \in Z_n, \textit{ with}\quad d_i\neq d_j \text{ for } i\neq j,
\]
and a symmetric matrix $Q\in V_n$, consider the condition
\begin{equation}\label{eq:isospectral}
\charpoly_{D+tQ}(x)=\charpoly_D(x)\qquad \text{for all }t\in K.
\end{equation}
So the affine line $D+KQ$ is required to lie in the fiber of the
morphism $\chi$  over the $K$-valued point $(d_1,\dots,d_n)\in \mathbb{A}^n$. Since $D$ has simple
spectrum, it is a regular semisimple point, but we restrict the direction
$Q$ to the linear space $V_n$. The resulting condition is a
homogeneous system of $\frac{n(n+1)}{2}$ polynomial equations in the $\frac{n(n+1)}{2}$ entries of $Q$ whose
coefficients depend algebraically on the diagonal entries of $D$.

This system occurs naturally in the geometric invariant theory of
multiple quadratic forms. The primary motivation for this paper comes from the recent work on the $\delta$-invariant theory of
Hecke correspondences by Buium and Vasiu where they  study the congruence action of
$\operatorname{SL}_g$ on tuples of symmetric matrices and the
$\Theta$-invariants obtained by polarizing the determinant
\cite[Section~4.2.6 and 4.2.9]{buiumvasiu2026}. The result obtained in this paper, together with the analysis in \cite[Section~4.2.12]{buiumvasiu2026} will allow us to compute the degree of a field extension that appears naturally in invariant theory.  This invariant-theoretic
application is our main motivation to study the isospectral condition over fields of positive characteristic.

Isospectral questions also have a classical analytic history.
Lax  exhibited a general mechanism for producing isospectral
deformations- an evolution satisfying $\dot L=[P,L]$ preserves the spectrum
of $L$ \cite[Section~1]{lax1968}. 

There is a  contrast between the real and algebraically closed
settings. The isospectral condition forces $Q$ to be nilpotent. Over
$\mathbb R$, a symmetric nilpotent matrix must vanish, so there are no
nontrivial deformations of this kind. Over an algebraically closed field,
however, nonzero symmetric nilpotent matrices do exist. Nilpotency alone
therefore does not solve the problem over algebraically closed fields. We must also use the  symmetry of $Q$, the simple spectrum of $D$, and the requirement
that the characteristic polynomial remains fixed for every $t\in K$.

The cases $n\in \{2,3\}$ are treated in \cite[Section~4.2.12]{buiumvasiu2026}; see also \cite[Proposition~5.1]{asokan2026}. In this paper we will prove some general results and discuss dimension four in detail. Our main
result shows that a generic diagonal matrix in dimension four admits no
nonzero symmetric direction for which the entire affine line is
isospectral over fields $K$ of characteristic $p\neq 2$ .

\medskip

The following is the dimension-four case of the main result of the
 author's dissertation \cite[Theorem~4.4 and Section~5.1.1]{asokan2026}.
\begin{theorem}\label{thm:main}
Let $K$ be an algebraically closed field with $\operatorname{char}(K)\neq 2$. For a generic diagonal matrix
\[
D=\diag(d_1,d_2,d_3,d_4)
\]
with pairwise distinct diagonal entries, the condition
\[
\charpoly_{D+tQ}(x)=\charpoly_D(x)\qquad \text{for all }t\in K
\]
for $Q\in V_4$ implies $Q=0$.
\end{theorem}

\medskip
\noindent
The proof begins with general consequences of the isospectral condition,
including the nilpotency of $Q$ and the vanishing of its diagonal entries. We then reduce the generic statement to finding one
diagonal matrix for which no nonzero symmetric perturbation exists. For a
convenient diagonal matrix in dimension four, we reduce the resulting
homogeneous equations to a smaller system and rule out its remaining
solutions by a finite-field computation. Finally, we use computer-assisted resultant calculations to
handle the exceptional positive characteristics and establish the theorem in
every characteristic different from $2$.\\
In Section~\ref{sec:invariants} we will combine Theorem~\ref{thm:main} with the degree criterion of Buium--Vasiu \cite[Theorem~4.2.22]{buiumvasiu2026} to compute the degree of a field extension appearing in the invariant theory of quadratic forms.
\medskip
This paper is based on the author's doctoral dissertation
\cite[Chapters~4-5]{asokan2026}. We present the dimension-four isospectral result
and its application to the invariant theory of triples of quadratic forms.\\
The behavior of isospectral deformations in dimensions \(n\geq 5\) differs significantly from that in dimension four. In particular, nontrivial symmetric isospectral deformations exist for generic diagonal matrices in higher dimensions. We will study these matrices and provide explicit constructions in an upcoming paper.

\section{ General properties of isospectral deformations}
This section presents the general reductions from
\cite[Section~4.2]{asokan2026}. We retain the arguments needed for the
dimension-four analysis and state the affine normalization in arbitrary
dimension.

\begin{definition}
For a diagonal matrix $D\in Z_n$, define
\[
S_{D,n}:=\left\{Q\in V_n:
\charpoly_{D+tQ}(x)=\charpoly_D(x)\text{ for all }t\in K\right\}.
\]
Equivalently, \[
S_{D,n}:=\left\{Q\in V_n:
\charpoly_{\lambda D+Q}(x)=\charpoly_{\lambda D}(x)\text{ for all }\lambda\in K\right\}.
\]
\end{definition}
We write $$\chi_{\lambda D+Q}(x)-\chi_{\lambda D}(x)=\sum_{a+b\leq n}P^{a,b}(D,Q)\lambda^ax^b,$$ where $a,b\in \mathbb{Z}_{\geq 0}$ and $P^{a,b}(D,Q)\in K[x_1,\dots,x_n,T_{ij} \mid 1\le i\le j\le n]$ is a homogeneous polynomial in $\frac{n(n+1)}{2}+n$ indeterminates. The polynomial $P^{a,b}(D,Q)$ is bihomogeneous of degree $a$ in the entries of $D$ and degree $n-a-b$ in the entries of $Q$. In particular, it is homogeneous of total degree $n-b$ in the entries of $D$ and $Q$. For a fixed $D=D_0$, we obtain the following system of $\frac{n(n+1)}{2}$ equations in as many variables:
\begin{equation}\label{main}
    P^{a,b}(D_0,Q)=0 \:\:\forall\:\: a,b\in \mathbb{Z}_{\geq 0}, a+b\leq n-1.
\end{equation}
By $I_D$ we denote the homogeneous  ideal of $K[T_{ij}|1\le i\le j\le n]$ generated by the set \\$\{ P^{a,b}(D,T_{11},\dots,T_{n,n}) \:\:\forall\:\: a,b\in \mathbb{Z}_{\geq 0}, a+b\leq n\}.$ We define the affine variety $\mathbf{S}_{D,n}:=\operatorname{Spec}(K[T_{ij}\mid 1\le i\le j\le n]/I_D)$ whose $K$-valued points are given by $S_{D,n}.$
\subsection{Symmetry with respect to permutations and translations}
Lemmas~\ref{Sym} and~\ref{5.2} will help simplify the system of equations (\ref{main}).
\begin{lemma}[\cite{asokan2026}, Proposition~4.7]\label{Sym}
Let $D\in Z_n$ and $Q\in V_n$ satisfy
\[
\chi_{D+tQ}(x)=\chi_D(x) \qquad \forall\,t\in K .
\]
For $\pi\in S_n$, define $\pi(A)$ by $(\pi(A))_{ij}=A_{\pi(i)\pi(j)}$.
Then
\[
\chi_{\pi(D)+t\pi(Q)}(x)=\chi_{\pi(D)}(x) \qquad \forall\,t\in K .
\]
\end{lemma}

\begin{proof}
Let $J_\pi$ be the permutation matrix corresponding to $\pi$. Then
\[
J_\pi(D+tQ)J_\pi^{-1}=\pi(D)+t\pi(Q).
\]
Since characteristic polynomials are invariant under similarity,
\[
\chi_{\pi(D)+t\pi(Q)}(x)=\chi_{D+tQ}(x).
\]
By hypothesis $\chi_{D+tQ}(x)=\chi_D(x)$, and since $\pi(D)$ is similar to $D$
we have $\chi_{\pi(D)}(x)=\chi_D(x)$. So
\[
\chi_{\pi(D)+t\pi(Q)}(x)=\chi_{\pi(D)}(x).
\]
\end{proof}
\begin{lemma}\label{TI}
    Let $D=\operatorname{diag}(d_1,\dots,d_n)$ be in $Z_n$ and $\lambda \in K.$ Let $D_{\lambda}=D+\lambda I_n.$ Then we have $S_{D,n}=S_{D_\lambda,n}.$ 
\end{lemma}
\begin{proof}
For any $Q\in V_n$ and $t\in K$, we have
$$
\chi_{D_\lambda+tQ}(x)
=\det(xI_n-D-\lambda I_n-tQ)\\
=\det((x-\lambda)I_n-D-tQ)\\
=\chi_{D+tQ}(x-\lambda).$$
Similarly,
\[
\chi_{D_\lambda}(x)=\chi_D(x-\lambda).
\]
So,
\[
\chi_{D_\lambda+tQ}(x)=\chi_{D_\lambda}(x)
\quad\text{for all }t\in K
\]
if and only if
\[
\chi_{D+tQ}(x-\lambda)=\chi_D(x-\lambda)
\quad\text{for all }t\in K.
\]
Since the substitution $x\mapsto x-\lambda$ is an
automorphism of $K[x]$, this is equivalent to
\[
\chi_{D+tQ}(x)=\chi_D(x)
\quad\text{for all }t\in K.
\]
Thus $Q\in S_{D_\lambda,n}$ if and only if
$Q\in S_{D,n}$, proving the lemma.
\end{proof}
The affine normalization below is the arbitrary-dimensional form of
\cite[Lemma~5.2]{asokan2026}; its translation step is given in
Lemma~\ref{TI}.
\begin{lemma}\label{5.2}
Let $K$ be algebraically closed and
$
D=\operatorname{Diag}(d_1,\ldots,d_n)\in Z_n
$ with distinct diagonal entries.
If $Q\in V_n$ satisfies
$
\chi_{D+tQ}(x)=\chi_D(x)\quad\text{for all }t\in K,
$
then, there exist
$
D'=\operatorname{Diag}(0,1,a_3,\ldots,a_n)\in Z_n,
$
with
$
a_i=\frac{d_i-d_1}{d_2-d_1},
$
such that
$$
\chi_{D'+tQ}(x)=\chi_{D'}(x)\quad\text{for all }t\in K.
$$
Moreover, the $a_i$ are pairwise distinct  and lie in $K\setminus{0,1}$.
\end{lemma}

\begin{proof}
We can replace $D$ by $D-d_1I_n$ using Lemma \ref{TI}. Let $c=d_2-d_1\neq0$ and set
$$
D'=c^{-1}(D-d_1I_n).
$$
Then
$$
\det(xI_n-D'-tQ)
=c^{-n}\det(cxI_n-(D-d_1I_n)-ctQ),
$$
which equals
$$
c^{-n}\det(cxI_n-(D-d_1I_n))
=\det(xI_n-D')
$$
since the original identity holds for all $t\in K$. The lemma follows since $d_i$ are distinct.
\end{proof}
\subsection{Nilpotence}

\begin{proposition}[\cite{asokan2026}, Proposition~4.5]\label{prop:nilpotent}
If $Q\in S_{D,n}$, then $Q$ is nilpotent.
\end{proposition}

\begin{proof}
This follows immediately by setting $\lambda=0$ in the definition of $S_{D,n}$ and applying the Cayley-Hamilton theorem. 
\end{proof}

\begin{remark}
Over $\mathbb{R}$, the only symmetric nilpotent matrix is the zero matrix, so the problem is only interesting over fields where nonzero symmetric nilpotent matrices can occur.
\end{remark}
\subsection{Vanishing of the diagonal entries of $Q\in S_{D,n}$}
The vanishing of the diagonal entries is established in
\cite[Proposition~4.9(1)]{asokan2026}. We give the coefficient argument
below; compare also \cite[Proposition~7.1]{asokan2026}.
\begin{proposition}
Let $D=\operatorname{Diag}(d_1,\ldots,d_n)\in Z_n$ have pairwise distinct diagonal entries, and let $Q=(q_{ij})\in S_{D,n}$. Then
$q_{ii}=0 \quad\text{for all }1\leq i\leq n.$
\end{proposition}
\begin{proof}
From
$
\chi_{D+tQ}(x)=\chi_D(x)$
we obtain
\[
\operatorname{Tr}\bigl((D+tQ)^k\bigr)=\operatorname{Tr}(D^k)
\qquad\text{for all }k\geq1.
\]
Taking the coefficient of $t$, we get
\[
\operatorname{Tr}
\bigl(QD^{k-1}+DQD^{k-2}+\cdots+D^{k-1}Q\bigr)=0.
\]
By cyclicity of trace, each term is equal to
$\operatorname{Tr}(QD^{k-1})$. So
$$
k\operatorname{Tr}(QD^{k-1})
=
k\sum_{i=1}^n q_{ii}d_i^{k-1}=0.
$$
If $\operatorname{char}(K)=0$, take $k=1,\ldots,n$. We obtain a Vandermonde system with determinant\\
$
\prod_{1\leq i<j\leq n}(d_j-d_i)\neq0.
$
Hence $q_{ii}=0$ for all $1\leq i\leq n$.\\
Suppose now that $\operatorname{char}(K)=p>0$. Take
$
k=1,1+p,\ldots,1+(n-1)p.
$
Since each of these integers is equal to $1$ modulo $p$, we obtain
$
\sum_{i=1}^n q_{ii}(d_i^p)^m=0
\quad\text{for }m=0,\ldots,n-1.
$
This is again a Vandermonde system, with determinant
$
\prod_{1\leq i<j\leq n}(d_j^p-d_i^p).
$ 
Since the $d_i$ are pairwise distinct, and the Frobenius endomorphism is injective this determinant is
nonzero. Therefore $q_{ii}=0$ for all $i$.
\end{proof}
\subsection{Reduction from generic $D$ to one special $D$}
\begin{proposition}[\cite{asokan2026}, Proposition~4.6]\label{prop:specialization}
Assume that there exists \(D_0\in \mathbf{Z}_n(K)\) such that there is no nonzero
\(Q\in V_n\) satisfying
\[
\chi_{D_0+tQ}(x)=\chi_{D_0}(x)\qquad \forall t\in K.
\]
 Then there exists a non-empty Zariski open subset
\(U\subset \mathbf{Z}_n\) such that for every \(D\in U(K)\), there is no nonzero
\(Q\in V_n\) satisfying
\[
\chi_{D+tQ}(x)=\chi_D(x)\qquad \forall t\in K.
\]
\end{proposition}

\begin{proof}
Let $
\mathbf{Y}_n:=\left\{(D,[Q])\in \mathbf{Z}_n\times \mathbb{P}\mathbf{V}_n \;\middle|\;
\chi_{D+tQ}(x)=\chi_D(x)\text{ for all }t\in K\right\}$ be the subvariety of $\mathbf{Z}_n\times \mathbb{P}\mathbf{V}_n$ defined by the system of homogeneous equations \ref{main}, and let
\[
\pi:\mathbf{Y}_n\to \mathbf{Z}_n,\qquad (D,[Q])\mapsto D
\]
be the projection morphism. Since \(\mathbf{Z}_n\times \mathbb{P}\mathbf{V}_n\to \mathbf{Z}_n\) is projective,
the morphism \(\pi\) is projective.

For any \(D\in \mathbf{Z}_n(K)\), the fiber \(\pi^{-1}(D)\) is the projective scheme of
nonzero solutions \(Q\in V_n\), up to scalar, of equation (\ref{eq:isospectral}),
$\chi_{D+tQ}(x)=\chi_D(x)\:\; \forall t\in K.$
In particular,

$\pi^{-1}(D)(K)=\varnothing$
if and only if there is no nonzero \(Q\in V_n\) satisfying the above polynomial equality.
By assumption, for \(D_0\in \mathbf{Z}_n(K)\) there is no nonzero \(Q\in V_n\) such that
\[
\chi_{D_0+tQ}(x)=\chi_{D_0}(x)\qquad \forall t\in K.
\]
Then,
$
\pi^{-1}(D_0)(K)=\varnothing.$
Since \(K\) is algebraically closed, this implies that \(\pi^{-1}(D_0)\) is empty as a scheme.
Hence $ D_0\notin \pi(\mathbf{Y}_n).$
Because \(\pi\) is projective, it is proper, and therefore \(\pi(\mathbf{Y}_n)\) is closed in
\(\mathbf{Z}_n\). It follows that
\[
U:=\mathbf{Z}_n\setminus \pi(\mathbf{Y}_n)
\]
is a non-empty Zariski open subset of \(\mathbf{Z}_n\). For every \(D\in U(K)\), the fiber
\(\pi^{-1}(D)\) is empty, so there is no nonzero \(Q\in V_n\) satisfying
\[
\chi_{D+tQ}(x)=\chi_D(x)\qquad \forall t\in K.
\]
This proves the proposition.
\end{proof}

\begin{remark}
This proposition allows the dimension-four problem to be reduced to a convenient explicit choice of diagonal entries.
\end{remark}

For $n\leq 3$, direct computation gives $S_{D,n}=\{0\}$; see \cite[Section 4.2.12]{buiumvasiu2026} and \cite[Proposition~5.1]{asokan2026}. We now turn to the case $n=4$.

\section{Dimension four}
We present the dimension-four argument of
\cite[Section~5.1.1]{asokan2026}. The proof proceeds through the
normalization of $D$, the non-vanishing of the off-diagonal entries,
and the elimination and resultant computations developed there.

Throughout this section, 
\[
D\in Z_4,
\qquad
Q=(q_{ij})\in V_4,
\]
where $D$ has pairwise distinct diagonal entries and $Q\in S_{D,4}$.
\begin{remark}
    Without loss of generality we may assume that the diagonal matrix $D$ is of the form $\diag(0,1,a,b).$
\end{remark}
For $Q\in S_{D,4}$, we obtain the following system of equations $C_4^{a,b}$,
corresponding to \cite[equations~(5.7)--(5.12)]{asokan2026}:
\begin{align}
q_{14}^2 q_{23}^2 + q_{13}^2 q_{24}^2 - 2q_{13}q_{14}q_{23}q_{24}
+ q_{12}^2 q_{34}^2 - 2q_{12}q_{14}q_{23}q_{34} - 2q_{12}q_{13}q_{24}q_{34} &= 0, \label{eq:C1'}\\
2q_{23}q_{24}q_{34} + 2q_{13}q_{14}q_{34} + 2q_{12}q_{14}q_{24} + 2q_{12}q_{13}q_{23} &= 0, \label{eq:C2'}\\
q_{13}q_{14}q_{34} + a q_{12}q_{14}q_{24} + b q_{12}q_{13}q_{23} &= 0, \label{eq:C3'}\\
a q_{14}^2 + b q_{13}^2 + ab q_{12}^2 &= 0, \label{eq:C4'}\\
(a+b)q_{12}^2 + (1+b)q_{13}^2 + (1+a)q_{14}^2 + b q_{23}^2 + a q_{24}^2 + q_{34}^2 &= 0, \label{eq:C5'}\\
q_{12}^2 + q_{13}^2 + q_{14}^2 + q_{23}^2 + q_{24}^2 + q_{34}^2 &= 0. \label{eq:C6'}
\end{align}

The next lemma uses the argument of \cite[Lemma~5.3]{asokan2026}.
\begin{lemma}\label{5.3}
Let the diagonal matrix,
\[
D=\diag(0,1,a,b),
\]
where $a,b\in K\setminus\{0,1\}$ and $a\neq b$ and let $Q=(q_{ij})\in S_{D,4}$ be nonzero. Then
\[
q_{ij}\neq 0 \qquad \text{for all } 1\le i<j\le 4.
\]
\end{lemma}

\begin{proof}
Suppose $q_{ij}=0$ for some $i\neq j$. By symmetry, we may assume $q_{12}=0$.

From the equations (\ref{eq:C1'}),(\ref{eq:C3'}) and (\ref{eq:C4'}) we obtain
\begin{align}
q_{14}q_{23}-q_{13}q_{24}&=0, \label{eq:L1}\\
q_{13}q_{14}q_{34}&=0, \label{eq:L2}\\
a q_{14}^2 + b q_{13}^2&=0. \label{eq:L3}
\end{align}

\textbf{Case 1: $q_{34}\neq 0$.}
Then \eqref{eq:L2} gives $q_{13}q_{14}=0$. Since $a,b\neq 0$, \eqref{eq:L3} implies $q_{13}=q_{14}=0$. So the lower-right $3\times 3$ principal block of $Q$ is nonzero and satisfies
\[
\chi_{\widetilde{D}+t\widetilde{Q}}(x)=\chi_{\widetilde{D}}(x),
\]
where $\widetilde{D}=\diag(1,a,b)$. This contradicts the $n=3$ case.

\textbf{Case 2: $q_{34}=0$.}
If $q_{13}q_{14}=0$, then again \eqref{eq:L3} gives $q_{13}=q_{14}=0$, and the same reduction to the $n=3$ case yields a contradiction.

So $q_{13}q_{14}\neq 0$. From \eqref{eq:L3} we obtain
\[
q_{13}^2 = -\frac{a}{b} q_{14}^2,
\]
and from \eqref{eq:L1},
\[
q_{23} = \frac{q_{13}q_{24}}{q_{14}}.
\]
From this we have $q_{23}^2=\frac{q_{13}^2q_{24}^2}{q_{14}^2}=\frac{-a}{b}q_{24}^2.$
Substituting into equation \eqref{eq:C5'} yields
\[
\left(1-\frac{a}{b}\right) q_{14}^2 = 0,
\]
a contradiction since $a\neq b$.

So all $q_{ij}\neq 0$.
\end{proof}
From Lemma \ref{5.3} we get the following proposition.
\begin{proposition}
    Let $K$ be an algebraically closed field with $\operatorname{char}(K)\neq 2$,  and let $D\in Z_4$  be a  diagonal matrix with distinct diagonal entries. Then $S_{D,4}$ is a union of finitely many lines through origin in $V_4$ or $S_{D,4}=\{0\}$. 
\end{proposition}
\begin{proof}
    From Lemmas \ref{Sym} and \ref{5.2} we can assume $D=\operatorname{diag}(0,1,a,b)$ for $(a,b)\in (K-\{0\})^2.$ 
    From lemma \ref{5.3} we have that the projectivization $\mathbb{P}S_{D,4}$ is a Zariski closed subset that lies in the affine chart $q_{12}\neq 0$ in $\mathbb{P}V_4.$ So $\mathbb{P}S_{D,4}$ is affine and proper over $K$ and hence finite over $K.$ 
\end{proof}
Lemma~\ref{5.3} allows us to simplify the system $C_4^{a,b}$.
The elimination below follows the proof of
\cite[Proposition~5.4]{asokan2026}.

Set
\[
x:=q_{12}, \qquad y:=q_{13}, \qquad z:=q_{23},
\]
and
\[
u:=x^2,\qquad v:=y^2,\qquad w:=z^2.
\]
Also define
\[
R:=q_{14}q_{24},\qquad S:=q_{14}q_{34},\qquad T:=q_{24}q_{34},
\]
and
\[
\ell_1:=q_{14}^2,\qquad \ell_2:=q_{24}^2,\qquad \ell_3:=q_{34}^2.
\]

We first solve for $\ell_1,\ell_2,\ell_3$ in terms of $u,v,w$, obtaining
\begin{align}
\ell_1 &= -bu-\frac{b}{a}v, \label{eq:l1'}\\
\ell_2 &= (b-1)u+\frac{1-b}{a-1}w, \label{eq:l2'}\\
\ell_3 &= \frac{b-a}{a}v+\frac{b-a}{a-1}w. \label{eq:l3'}
\end{align}

Next, equations \eqref{eq:C1'}, \eqref{eq:C2'}, and \eqref{eq:C3'} give relations among $R,S,T$. Equation \eqref{eq:C1'} becomes
\[
2\sqrt{vw}\,R + 2\sqrt{uw}\,S + 2\sqrt{uv}\,T = P,
\]
where
\[
P:=w\ell_1+v\ell_2+u\ell_3.
\]
Using \eqref{eq:l1'}, \eqref{eq:l2'}, and \eqref{eq:l3'}, we compute
\begin{equation}\label{P}P=\frac{ab-2a+b}{a}uv+\frac{-ab-a+2b}{a-1}uw+\frac{-2ab+a+b}{a(a-1)}vw.
\end{equation}
Equation \eqref{eq:C2'} becomes
\[
\sqrt{u}\,R+\sqrt{v}\,S+\sqrt{w}\,T=-\sqrt{uvw},
\]
and \eqref{eq:C3'} becomes
\[
a\sqrt{u}\,R+\sqrt{v}\,S=-b\sqrt{uvw}.
\]
As a result,
\[
\begin{pmatrix}
2\sqrt{vw} & 2\sqrt{uw} & 2\sqrt{uv}\\
\sqrt{u} & \sqrt{v} & \sqrt{w}\\
a\sqrt{u} & \sqrt{v} & 0
\end{pmatrix}
\begin{pmatrix}
R\\
S\\
T
\end{pmatrix}
=
\begin{pmatrix}
P\\
-\sqrt{uvw}\\
-b\sqrt{uvw}
\end{pmatrix}.
\]
Let
\[
M:=
\begin{pmatrix}
2\sqrt{vw} & 2\sqrt{uw} & 2\sqrt{uv}\\
\sqrt{u} & \sqrt{v} & \sqrt{w}\\
a\sqrt{u} & \sqrt{v} & 0
\end{pmatrix},
\qquad
B:=\det(M)=(2-2a)uv+2auw-2vw.
\]
A direct computation gives
\[
\operatorname{Adj}(M)=
\begin{pmatrix}
-\sqrt{vw} & 2\sqrt{uv} & -2\sqrt{u}(v-w)\\
a\sqrt{uw} & -2a\sqrt{v}u & 2\sqrt{v}(u-w)\\
\sqrt{uv}(1-a) & 2\sqrt{w}(au-v) & 2\sqrt{w}(v-u)
\end{pmatrix}.
\]
So
\[
B
\begin{pmatrix}
R\\
S\\
T
\end{pmatrix}
=
\operatorname{Adj}(M)
\begin{pmatrix}
P\\
-\sqrt{uvw}\\
-b\sqrt{uvw}
\end{pmatrix},
\]
and therefore
\begin{align}
BR &= \sqrt{vw}\Big((2b-2)uv-2buw-P\Big), \label{eq:BR'}\\
BS &= \sqrt{uw}\Big(aP+(2a-2b)uv+2bvw\Big), \label{eq:BS'}\\
BT &= \sqrt{uv}\Big((1-a)P+(2b-2a)uw+(2-2b)vw\Big). \label{eq:BT'}
\end{align}
Since
\[
R^2=\ell_1\ell_2,\qquad S^2=\ell_1\ell_3,\qquad T^2=\ell_2\ell_3,
\]
Equations \eqref{eq:BR'}, \eqref{eq:BS'}, and \eqref{eq:BT'} when squared gives
\begin{align}
vw\Big((2b-2)uv-2buw-P\Big)^2 &= B^2\Big(-bu-\frac{b}{a}v\Big)\Big((b-1)u+\frac{1-b}{a-1}w\Big), \label{eq:E1'}\\
uw\Big(aP+(2a-2b)uv+2bvw\Big)^2 &= B^2\Big(-bu-\frac{b}{a}v\Big)\Big(\frac{b-a}{a}v+\frac{b-a}{a-1}w\Big), \label{eq:E2'}\\
uv\Big((1-a)P+(2b-2a)uw+(2-2b)vw\Big)^2 &= B^2\Big((b-1)u+\frac{1-b}{a-1}w\Big)\Big(\frac{b-a}{a}v+\frac{b-a}{a-1}w\Big). \label{eq:E3'}
\end{align}

Substituting for $P$ from Equation \ref{P} we obtain
\begin{align}
vw\left(\frac{ab-b}{a}uv+\frac{-ab+a}{a-1}uw+\frac{2ab-a-b}{a(a-1)}vw\right)^2
&=
B^2\left(-bu-\frac{b}{a}v\right)\left((b-1)u+\frac{1-b}{a-1}w\right), \label{eq:F1'}\\
uw\left((ab-b)uv+\frac{-a^2b-a^2+2ab}{a-1}uw+\frac{a-b}{a-1}vw\right)^2
&=
B^2\left(-bu-\frac{b}{a}v\right)\left(\frac{b-a}{a}v+\frac{b-a}{a-1}w\right), \label{eq:F2'}\\
uv\left(\frac{-a^2b+2a^2-2a+b}{a}uv+(ab-a)uw+\frac{a-b}{a}vw\right)^2
&=
B^2\left((b-1)u+\frac{1-b}{a-1}w\right)\left(\frac{b-a}{a}v+\frac{b-a}{a-1}w\right). \label{eq:F3'}
\end{align}
Since $Q\neq 0$, Lemma \ref{5.3} gives
\[
uvw\neq 0.
\]
We may therefore set
\[
U:=\frac{1}{u},\qquad V:=\frac{1}{v},\qquad W:=\frac{1}{w},
\]
and
\[
\widetilde{B}:=\frac{B}{uvw}=(2-2a)W+2aV-2U.
\]
Dividing \eqref{eq:F1'}, \eqref{eq:F2'}, and \eqref{eq:F3'} by $(uvw)^2(uv)(uw)$, $(uvw)^2(uv)(vw)$, and $(uvw)^2(uw)(vw)$ respectively, we obtain the system
\begin{align}
U^2\left(\frac{ab-b}{a}W+\frac{-ab+a}{a-1}V+\frac{2ab-a-b}{a(a-1)}U\right)^2
&-\widetilde{B}^2\left(-\frac{b}{a}U-bV\right)\left((b-1)W+\frac{1-b}{a-1}U\right)=0, \label{eq:L29'}\\
V^2\left((ab-b)W+\frac{-a^2b-a^2+2ab}{a-1}V+\frac{a-b}{a-1}U\right)^2
&-\widetilde{B}^2\left(-\frac{b}{a}U-bV\right)\left(\frac{b-a}{a}W+\frac{b-a}{a-1}V\right)=0, \label{eq:L30'}\\
W^2\left(\frac{-a^2b+2a^2-2a+b}{a}W+(ab-a)V+\frac{a-b}{a}U\right)^2
&-\widetilde{B}^2\left((b-1)W+\frac{1-b}{a-1}U\right)\left(\frac{b-a}{a}W+\frac{b-a}{a-1}V\right)=0. \label{eq:L31'}
\end{align}
We denote this reduced system by $L(a,b)$. Every solution of the system $C_4^{a,b}$ will induce a solution for system $L(a,b).$ While the systems are not equivalent, to show that $C_4^{a,b}$ has no nonzero solutions, it suffices to show $L(a,b)$ has no nonzero solutions.\\
 We choose $(a,b)=(3,4)$. Then $L(3,4)$ is the following system,
\begin{align}
U^2\left(\frac{8}{3}W-\frac{9}{2}V+\frac{17}{6}U\right)^2
-( -4W+6V-2U )^2\left(-\frac{4}{3}U-4V\right)\left(3W-\frac{3}{2}U\right)&=0, \label{eq:G1'}\\
V^2\left(8W-\frac{21}{2}V-\frac{1}{2}U\right)^2
-( -4W+6V-2U )^2\left(-\frac{4}{3}U-4V\right)\left(\frac{1}{3}W+\frac{1}{2}V\right)&=0, \label{eq:G2'}\\
W^2\left(-\frac{20}{3}W+9V-\frac{1}{3}U\right)^2
-( -4W+6V-2U )^2\left(\frac{1}{3}W+\frac{1}{2}V\right)\left(3W-\frac{3}{2}U\right)&=0. \label{eq:G3'}
\end{align}
The finite-field computation below is the certificate from
\cite[equations~(5.46)--(5.51) and Remark~5.5]{asokan2026}.
\begin{lemma}
    \label{P2}
    The homogeneous system of equations $L(3,4)$ has no solutions in $\mathbb{P}^2_{\bar{\mathbb{F}}_{11}}=\operatorname{Proj}\bar{\mathbb{F}}_{11}[U,V,W]$ outside the hypersurface defined by $\widetilde{B}=0.$
\end{lemma}
\begin{proof}
    We restrict to the affine chart $\widetilde{B}=1$. Reducing the resulting system modulo $11$, we obtain the following three equations in $\mathbb{F}_{11}[V,W]$:
\begin{align}
f_1 &:= V^4+10V^3W+4V^3+4V^2W^2+3V^2W+7V^2
\notag\\
&\qquad
+5VW^3+5VW^2+VW+10V
+3W^4+8W^3+7W^2+9W+3, \label{eq:f1'}
\\[4pt]
g_1 &:= V^4+4V^3W+5V^3+4V^2W^2+10V^2W+2V^2
\notag\\
&\qquad
+5VW+7V+4W^2+W,\label{eq:g1'}
\\[4pt]
h_1 &:= 9V^2W^2+5V^2+3VW^3+10VW^2+4VW+V
\notag\\
&\qquad
+3W^4+9W^3+2W^2+8W\label{eq:h1'}.
\end{align}
We will show that $f_1$, $g_1$, and $h_1$ have no common zero. Consider
\begin{align}
M_1 &:= 5V^3W-V^2W^2+4W^4+4V^3+V^2W-5W^3
-4V^2+W^2+2V-3W+4,
\\[4pt]
M_2 &:= -5V^3W+4V^2W^2+5VW^3-4V^3+2V^2W+VW^2
 +3W^3-3V^2-3VW+W^2-2V+5W+5,
\\[4pt]
M_3 &:= 3V^2W^2+VW^3-4W^4+2VW^2-W^3
 +3V^2+3VW+5W^2-V-4W-4.
\end{align}
A direct computation gives
\[
f_1M_1+g_1M_2+h_1M_3=1
\]
in $\mathbb{F}_{11}[V,W]$. Hence $f_1$, $g_1$, and $h_1$ have no common zero. This proves the lemma.\\
The polynomials $M_1,M_2$ and $M_3$ were obtained by using the computer algebra system Macaulay2 \cite{Macaulay2}.
\end{proof}
 Now we have the following proposition in fields $K$ of characteristic $11.$ 
 \begin{proposition}\label{Ch11}
Suppose that $\operatorname{char}(K)$=11, and let
$D=\operatorname{Diag}(0,1,3,4).$
If $Q\in S_{D,4}$
then $Q=0$.
\end{proposition} 
\begin{proof} From Lemma \ref{P2} we rule out possible solutions away from the hypersurface. 
 So now we consider the hypersurface defined by $\widetilde{B}=0$ in $\mathbb{P}^2_{\bar{\mathbb{F}}_{11}}.$  The system of equations $L(3,4)$ reduces to a linear system with the unique projective solution $[U:V:W]=[1:3:4].$ From this point we get the possible values of $q_{12},\dots,q_{34}$ as \[
q_{12}=\pm1,\quad q_{13}=\pm2,\quad q_{23}=\pm5,\quad
q_{14}=\pm3,\quad q_{24}=\pm2,\quad q_{34}=\pm1.
\]
 For every choice of signs, the left-hand side of Equation (\ref{eq:C3'}) is
of the form $\pm6\pm7\pm7$ and is not  zero in $\mathbb F_{11}.$ So we do not get a nonzero solution for $C_4^{3,4}.$

\end{proof}
Now we are ready to prove the theorem \ref{thm:main}.
\subsection{Proof of  Theorem \ref{thm:main}}
We complete the proof with the resultant computations in
\cite[Section~5.1.1, following Definition~5.6]{asokan2026}.
We recall the following result from elimination theory.
\begin{proposition} Let $A$ be a domain and $L$ be an algebraic closure of the field of fractions of $A.$
    Let  $f_0,\dots,f_n$ be homogeneous polynomials in $A[x_0,\dots,x_n]$ with the Macaulay's resultant  $\mathcal{R}(f_0,\dots,f_n)\in A.$ Then $\mathcal{R}(f_0,\dots,f_n)=0$ if and only if \[\operatorname{Proj}\!\bigl(L[x_0,\dots,x_n]/(I)\bigr) \neq \varnothing.
\] 
\end{proposition}
\begin{proof}
    This is the classical vanishing criterion for the homogeneous resultant; see \cite{vanderWaerden1950} and the formulation in \cite[Section~2.3.1, Theorem~2.32]{asokan2026}.
\end{proof}
\begin{definition}
    Consider the homogeneous system $C_4^{a,b}$. We denote by $R_p^{a,b}$ the Macaulay resultant of $C_4^{a,b}$ modulo $p$, and by $R^{a,b}:=R_0^{a,b}$ the Macaulay resultant of $C_4^{a,b}$ in characteristic zero.
\end{definition}
From Proposition \ref{Ch11}  we have that $R_{11}^{3,4}\neq 0$ and hence $R^{3,4}\in \mathbb{Z}$ is not zero. Let $R^{3,4}=\prod_{i=1}^r p_i^{e_i}$ where $p_i$ are primes in $\mathbb{Z}.$ Then for all $p\notin \{p_1,\dots,p_r\}$ we have $R_p^{3,4}\neq 0.$ From Proposition \ref{prop:specialization}, for a generic $D$ in the open subset of $Z_4$ we have $S_{D,4}=\{0\}.$\\
Using the computer algebra system Macaulay2 \cite{Macaulay2}, we compute the Macaulay resultant of $C_4^{3,4}$ and obtain $R^{3,4}=2^{296}\times3^{112}\times17^{32}$. Similarly, we obtain $R^{4,5}=2^{352}\times3^{56}\times5^{56}\times31^{32}$. This implies that for fields $K$ of characteristic $p\notin \{2,3\}$ we have $S_{D,4}=\{0\}$  for a generic $D.$ 

In characteristic $3$, consider the field $\mathbb{F}_3[y]/(y^2+1)\cong \mathbb{F}_9$. For $D=\operatorname{Diag}(0,1,2,y)$, we obtain the following homogeneous system:
\[
C_4^{2,y}:
\begin{cases}
q_{14}^2q_{23}^2 + q_{13}^2q_{24}^2 - 2q_{13}q_{14}q_{23}q_{24}
+ q_{12}^2q_{34}^2 - 2q_{12}q_{14}q_{23}q_{34} - 2q_{12}q_{13}q_{24}q_{34} = 0, \\[6pt]

2q_{23}q_{24}q_{34} + 2q_{13}q_{14}q_{34}
+ 2q_{12}q_{14}q_{24} + 2q_{12}q_{13}q_{23} = 0, \\[6pt]

q_{13}q_{14}q_{34} + 2\,q_{12}q_{14}q_{24} + y\,q_{12}q_{13}q_{23} = 0, \\[6pt]

2\,q_{14}^2 + y\,q_{13}^2 + 2y\,q_{12}^2 = 0, \\[6pt]

(2+y)q_{12}^2 + (1+y)q_{13}^2 + 3q_{14}^2
+ y\,q_{23}^2 + 2\,q_{24}^2 + q_{34}^2 = 0, \\[6pt]

q_{12}^2 + q_{13}^2 + q_{14}^2
+ q_{23}^2 + q_{24}^2 + q_{34}^2 = 0.
\end{cases}
\]
The Macaulay resultant of the system $C_4^{2,y}$ is $1$. Hence $S_{D,4}=\{0\}$ for generic $D$ over every field of characteristic $p\neq 2$. This completes the proof of the Theorem\ref{thm:main}.
\subsection{On the generic open set}
Let
$
\Omega_p=
\left\{(a,b)\in\mathbb A^2_K:
a,b\notin\{0,1\},\ a\neq b,\ R_p^{a,b}\neq0
\right\}.
$
On the principal open subset $\mathcal{O}(d_2-d_1)\subseteq Z_4$, consider
the morphism
\[
\nu_{12}(D)=
\left(
\frac{d_3-d_1}{d_2-d_1},
\frac{d_4-d_1}{d_2-d_1}
\right).
\]
The corresponding open subset of $Z_4$ is
\[
\mathcal U_p:=\nu_{12}^{-1}(\Omega_p).
\]
For all $D\in \mathcal{U}_p$ we have $S_{D,4}=\{0\}.$
\section{Application to invariant theory}\label{sec:invariants}
We now apply Theorem~\ref{thm:main} to the invariant theory of triples
of quadratic forms. Using the construction of Buium-Vasiu
\cite[Sections~4.2.11--4.2.12]{buiumvasiu2026}, we compute the
degree of the extension from the field generated by determinant
polarizations to the full invariant field.
Throughout this section, $K$ is algebraically closed of characteristic
$p\neq2$.

\subsection{Invariant rings and determinant polarizations}
The group $\mathbf{SL}_n$ acts on $\mathbf V_n^{r+1}$ by simultaneous
congruence,
\[
\Lambda\cdot(Q_0,\ldots,Q_r)
=(\Lambda Q_0\Lambda^t,\ldots,\Lambda Q_r\Lambda^t).
\]
Write
\[
H_n^r(K)=K[\mathbf V_n^{r+1}]^{\mathbf{SL}_n},\qquad
X_n^{r+1}=\mathbf V_n^{r+1}//\mathbf{SL}_n
=\operatorname{Spec}(H_n^r(K)).
\]
The quotient construction is recalled in \cite[Section~2.2]{asokan2026}; see also \cite[Chapter~1]{MumfordFogartyKirwan1994}.
Our $H_n^r(K)$ is the total invariant algebra denoted
$H_{n,\mathrm{tot},K}^r$ in \cite{asokan2026} and
\cite{buiumvasiu2026}; it should not be confused with their
integer-graded subalgebra $H_{n,K}^r$.

Let
\[
\Delta(n,r)=\{(m_0,\ldots,m_r)\in\mathbb Z_{\geq0}^{r+1}:
 m_0+\cdots+m_r=n\}
\]
and define the $\Theta$-invariants by
\[
\det(y_0Q_0+\cdots+y_rQ_r)
=\sum_{\mathbf m\in\Delta(n,r)}
 \Theta_{\mathbf m}(Q_0,\ldots,Q_r)y_0^{m_0}\cdots y_r^{m_r}.
\]
Each coefficient is invariant under simultaneous congruence by
$\mathbf{SL}_n$, since $\det(\Lambda)^2=1$. Denote the subalgebra
generated by these coefficients by $H_n^{r,\Theta}(K)$. Thus
\[
H_n^{r,\Theta}(K)\subseteq H_n^r(K).
\]
For triples, the resulting field extension is finite, and we write
\begin{equation}\label{eq:d-invariant}
d_{2,n,p}:=
[\operatorname{Frac}(H_n^2(K)):
 \operatorname{Frac}(H_n^{2,\Theta}(K))].
\end{equation}
This degree measures the extension of invariant function fields determined by the $\Theta$-subalgebra.
For pairs, Buium--Vasiu prove that the total invariant algebra is the
polynomial algebra generated by the $n+1$ determinant coefficients
\cite[Theorem~4.2.8]{buiumvasiu2026}; see also
\cite[Theorem~3.8]{asokan2026}. In particular,
\[
X_n^2=\operatorname{Spec}(H_n^{1,\Theta}(K)).
\]

\subsection{The morphisms $\omega_n^{3,2}$ and $\psi_n^{3,2}$}
The projection onto the first two matrices induces a morphism
$X_n^3\to X_n^2$. Together with the quotient map on triples, it defines
\[
\omega_n^{3,2}:\mathbf V_n^3\longrightarrow
Y_n:=\mathbf V_n^2\times_{X_n^2}X_n^3.
\]
The inclusions of the determinant subalgebras similarly define
\[
\psi_n^{3,2}:\mathbf V_n^3\longrightarrow
Z_n^{\Theta}:=\mathbf V_n^2\times_{X_n^2}
\operatorname{Spec}(H_n^{2,\Theta}(K)).
\]
Here $Z_n^{\Theta}$ is a fiber product, not the space $\mathbf Z_n$
of diagonal matrices. The inclusion
$H_n^{2,\Theta}(K)\subseteq H_n^2(K)$ induces
$\varphi_n^{3,2}:Y_n\to Z_n^{\Theta}$, and the construction gives
\begin{equation}\label{eq:psi-factorization}
\psi_n^{3,2}=\varphi_n^{3,2}\circ\omega_n^{3,2}.
\end{equation}
This factorization relates the invariant-field degree to the degree of $\psi_n^{3,2}$; see also \cite[Section~4.1]{asokan2026}.

Let $\epsilon_n\in\{0,1\}$ be the remainder of $n$ modulo $2$.
Buium--Vasiu's generic-degree computation
\cite[Lemma~4.2.16]{buiumvasiu2026}, recalled in
\cite[Theorem~4.1]{asokan2026}, gives
\[
\deg(\omega_n^{3,2})=2^{n-2+\epsilon_n}.
\]
Moreover, $\omega_n^{3,2}$ is generically \'{e}tale when $p\neq2$.
As explained in \cite[Remark~4.2]{asokan2026}, the base change defining
$\varphi_n^{3,2}$ preserves the relevant function-field degree, so
\[
\deg(\varphi_n^{3,2})=d_{2,n,p}.
\]
Consequently, if $b_{p,n}:=\deg(\psi_n^{3,2})$, multiplicativity of
degrees in \eqref{eq:psi-factorization} yields
\begin{equation}\label{eq:b-versus-d}
b_{p,n}=2^{n-2+\epsilon_n}d_{2,n,p}.
\end{equation}
It remains to compute $b_{p,n}$.

\subsection{Reduction to the isospectral equations}
We now express the fibers of $\psi_n^{3,2}$ in terms of the isospectral equations; this reduction also appears in the dissertation
\cite[Section~4.1, equations~(4.1)--(4.5)]{asokan2026}.
Fix the first two matrices to be $I_n$ and a diagonal matrix $D$.
The fiber of $\psi_n^{3,2}$ through $(I_n,D,0)$ consists of the triples
$(I_n,D,Q)$ for which
\[
\Theta_{\mathbf m}(I_n,D,Q)=\Theta_{\mathbf m}(I_n,D,0)
\qquad(\mathbf m\in\Delta(n,2)).
\]
By the definition of the determinant polarizations, this is equivalent to
\[
\det(y_0I_n+y_1D+y_2Q)=\det(y_0I_n+y_1D).
\]
This condition is equivalent to
\[
\chi_{\lambda D+Q}(x)=\chi_{\lambda D}(x)
\quad\text{for all }\lambda\in K,
\]
or equivalently $Q\in S_{D,n}$. So the isospectral equations studied
in this paper describe the zero fiber of the coefficient map occurring
in the invariant-theoretic construction.

The degree criterion of Buium and Vasiu
\cite[Theorem~4.2.22]{buiumvasiu2026}, recalled in
\cite[Theorem~4.3]{asokan2026}, connects the two problems.
Assuming the invariant-field extension is separable, if $S_{D,n}=\{0\}$ for generic diagonal $D$, then
\begin{equation}\label{eq:b-degree-product}
b_{p,n}=\prod_{i=1}^n i^{\,n+1-i}.
\end{equation}

\subsection{The degree in dimension four}
Combining Theorem~\ref{thm:main} with the preceding degree formulas gives the following invariant-theoretic application
\cite[Section~5.2, Theorem~5.7]{asokan2026}.
\begin{theorem}\label{thm:invariant-degree}
Let $K$ be an algebraically closed field of characteristic $p\neq2$.
Assume that
\[
\operatorname{Frac}(H_4^{2,\Theta}(K))
\subseteq\operatorname{Frac}(H_4^2(K))
\]
is separable. Then
\[
[\operatorname{Frac}(H_4^2(K)):
 \operatorname{Frac}(H_4^{2,\Theta}(K))]=72.
\]
\end{theorem}
\begin{proof}
By Theorem~\ref{thm:main}, $S_{D,4}=\{0\}$ for generic diagonal $D$.
Applying \eqref{eq:b-degree-product} gives
\[
b_{p,4}=\prod_{i=1}^4 i^{5-i}
=1^4\cdot2^3\cdot3^2\cdot4=288.
\]
Since $\epsilon_4=0$, the degree of $\omega_4^{3,2}$ is
$2^{4-2}=4$. Equation~\eqref{eq:b-versus-d} therefore gives
\[
d_{2,4,p}=\frac{b_{p,4}}{\deg(\omega_4^{3,2})}
=\frac{288}{4}=72.
\]
\end{proof}

\begin{remark}
The conclusion holds in
characteristic zero and in sufficiently large positive characteristic,
by the separability result of Buium--Vasiu
\cite[Section~4.2.12]{buiumvasiu2026}. Although Theorem~\ref{thm:main}
holds in every characteristic different from $2$, extending the
field-degree conclusion beyond the separability hypothesis
requires a separate justification.
\end{remark}
\section*{Acknowledgements}

I would like to thank Prof. Adrian Vasiu for the careful reading of the initial draft and valuable suggestions.
\bibliographystyle{amsplain}
\bibliography{references}

\end{document}